\documentclass[11pt]{article}

\usepackage[letterpaper]{geometry}
\usepackage{microtype}
\usepackage{times}
\usepackage{graphicx}
\usepackage[colorlinks=true,allcolors=blue]{hyperref}

\title{Is Mathematics Obsolete?}

\author{Jeremy Avigad}

\date{February 14, 2025}

\begin{document}

\maketitle

\section{Mathematical Reasoning and Scientific Reasoning}

Displayed prominently in the Vatican, Raphael's fresco, \emph{The School of Athens}, portrays several philosophers and thinkers of ancient Greece, arrayed around Plato and Aristotle at the center. Plato, whose philosophy of knowledge is based on \emph{forms} that inhabit an extraterrestrial realm, points to the heavens. Aristotle, who took mathematical and scientific abstractions to be rooted in worldly experience, holds his hand out, open, palm facing downward. With Plato drawing our attention to the pristine beauty of mathematical forms, I imagine Aristotle saying ``Whoa there, Romeo, get your head out of the clouds. We have work to do, here on Planet Earth.''

\begin{figure}[ht!]
  \centering
  \includegraphics[width=70mm]{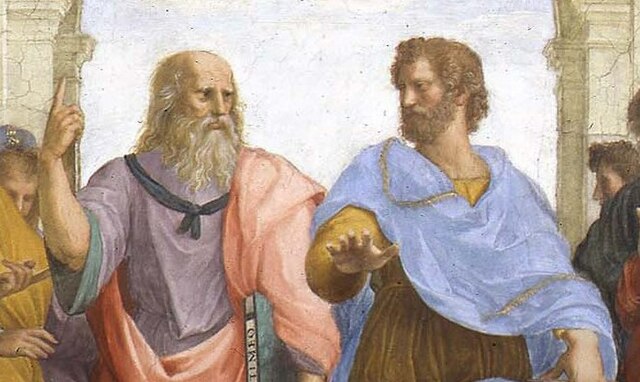}
  \caption{Detail from Raphael's \emph{School of Athens}.}
\end{figure}

The scene sums up the complementary appeals of mathematical and scientific thought. Mathematics deals with abstract objects, like numbers, that are not located in time or space, while science is about concrete objects, located in the physical world. Mathematics is supremely rational, relying on deductive reasoning to pass from axioms and hypotheses to conclusions, whereas science is fundamentally empirical, relying on inductive reasoning to pass from data to general laws. Mathematical claims are exact, while scientific claims are approximate. Mathematical deduction yields absolute certainty, while scientific induction yields, at best, hypotheses that are supported by the data but can be overturned by future evidence.

These contrasts have fascinated philosophers through the ages. Historians classify the early modern philosophers of the seventeenth and eighteenth centuries as \emph{rationalists} or \emph{empiricists}; the rationalists took knowledge to be grounded in abstract ideas, obtained by reason and independent of experience, whereas the empiricists held that the physical world is the ultimate source of all our knowledge. These characterizations evolved over the years. In the mid-twentieth century, the logical empiricist movement classified mathematical knowledge as \emph{analytic}, which is to say, true by virtue of the proper use of language and our shared linguistic framework, in contrast to the \emph{synthetic} knowledge of science.

What is at stake in the philosophical debates is often murky. Sometimes, it is a matter of \emph{ontology}, an account of the objects we talk about and how we should talk about them. Sometimes, it is a matter of \emph{epistemology}, sorting out the appropriate ways of justifying claims to knowledge. The philosophical positions, however, often come across as expressions of subjective preference. If you like mathematics, you are prone to view abstract mathematical knowledge as the most pristine form of knowledge and facts about the world as nothing more than pale shadows of ideal laws and mathematical truths. If you are more scientifically inclined, experience and observation are what counts, and you are more inclined to view mathematics as nothing more than the language we use to describe them.

Nobody can deny that mathematics and science need each other. One of the reasons we care about mathematics is that it gives us such a powerful means of thinking about the world, and however you feel about mathematics, it's impossible to imagine what contemporary science would look like without it. In the words of Immanuel Kant, thoughts without content are empty, and intuitions without concepts are blind. The question as to whether ideas or data come first is a chicken-and-egg problem: we build our conceptual scaffolding to make sense of our experiences, and that scaffolding, in turn, determines what can we make of them.

\section{Symbolic and Neural AI}

Given the excitement about neural networks and generative AI, I like to remind people that these represent only one of the two main approaches to artificial intelligence, nicely described in an \href{https://cacm.acm.org/research/deep-learning-for-ai/}{article} by Turing Award winners Yoshua Bengio, Yann LeCun, and Geoffrey Hinton:
\begin{quote}
 There are two quite different paradigms for AI. Put simply, the logic-inspired paradigm views sequential reasoning as the essence of intelligence and aims to implement reasoning in computers using hand-designed rules of inference that operate on hand-designed symbolic expressions that formalize knowledge. The brain-inspired paradigm views learning representations from data as the essence of intelligence and aims to implement learning by hand-designing or evolving rules for modifying the connection strengths in simulated networks of artificial neurons.
\end{quote}
Historically, the logic-inspired paradigm was there first. Early in the twentieth century, logicians wrote down axioms and rules for mathematical reasoning, and they developed algorithms to decide the truth of mathematical statements and search for mathematical proofs. They also determined some of the fundamental theoretical limits of symbolic reasoning, demonstrating that consistent proof systems are necessarily incomplete and that there are classes of mathematical problems with no algorithmic solution. These results were well established by the early 1940s, as the first digital computers were being built. A workshop in 1956, the \emph{Dartmouth Summer Research Project on Artificial Intelligence}, launched AI as a new field of research, and by the end of the decade, a number of practitioners had implemented logic-based reasoning algorithms.

Enthusiasm for symbolic AI ebbed and flowed in the subsequent decades. Rampant optimism in the 1960s gave way to the disappointing realization in the 1970s that simulating intelligent behavior was harder than AI's early proponents had hoped. The Japanese government's \emph{Fifth Generation Computer Systems} project, which aimed to develop massively parallel logic-based platforms as a basis for expert systems, brought new optimism in the early 1980s but gave way to renewed skepticism by the end of the decade. The 1990s thus inaugurated an ``AI winter'' in the US, with reduced government funding and support. Symbolic AI achieved landmark success in 1997, however, when IBM's Deep Blue beat the reigning world chess champion, Garry Kasparov, in a six-game tournament.

Despite these ups and downs, in the latter half of the twentieth century, the symbolic worldview reigned as the dominant paradigm in our understanding of thinking and reasoning. Cognitive science was a matter of explaining intelligent human behavior by sussing out the symbolic representations and rules that give rise to it. Linguistics was a matter of working out the symbolic grammars and rules that govern language generation and interpretation. These perspectives guided efforts in mechanization. Automated reasoning meant searching for proofs in logical calculi; designing an expert system meant encoding knowledge in rules and using them to draw inferences; natural language processing meant writing grammars and parsers to analyze and generate sentences; implementing systems of computer vision meant writing algorithms to detect corners, edges, and faces, and using that data to construct symbolic representations of the underlying objects and the relationships between them. In short, the common understanding of intelligence was that it was a matter of processing inputs into logic-based representations and operating on those. This is not to deny that the second half of the twentieth century also saw important advances in machine learning and neural networks; both were discussed at the Dartmouth workshop and were active research areas in the decades that followed. But despite the vicissitudes of the symbolic approach, the scales always seemed to tip in its favor.

Toward the end of the twentieth century, researchers in artificial intelligence retrenched and avoided grandiose claims about general intelligence. Practitioners focused on more specialized tasks like image recognition, natural language processing, and hardware and software verification. In this environment, machine learning, which aims to learn from data rather than hand-crafted rules, began to gain traction. Using empirical and statistical methods rather than logical ones, the field addressed tasks like classifying inputs and making predictions based on regularities in data.

The balance between the two approaches shifted decisively with highly publicized events like the success of IBM's Watson on \emph{Jeopardy!} in 2011, AlphaGo's defeat of Lee Sudol in 2016, and the arrival of the large language model GPT-3 in 2020. All of a sudden, AI was no longer about symbolic representations and reasoning but about big data and algorithms for learning. The goal was not to generate precise, rule-based justification and explanations but to use past experience and data to predict the most likely answers and promising strategies. Concepts were no longer localized but distributed opaquely across billions of parameters in a neural network. Rules of inference were no longer embedded into the design of a system but, instead, expected to emerge organically from the data and the training process. In the age-old dialogue between Plato and Aristotle, Aristotle had gained the upper hand.

\section{Is Mathematics Obsolete?}

To researchers like me, who had devoted their lives to the study of logical systems and symbolic methods, the sudden rise of machine learning and neural networks felt like an existential threat. Over and over again, we saw them come to dominate fields where symbolic methods once held sway. In cognitive science, symbolic models of human cognition gave way to neural models. The fifth-generation expert systems of the 1980s gave way to big data, foundation models, and neural networks. Deep Blue was eclipsed by AlphaZero, and symbolic approaches to language processing gave way to LLMs. Symbolic algorithms for image processing were crushed by neural networks.

This sea change was poignant at Carnegie Mellon, where I work. In the second half of the twentieth century, Herbert Simon, a Nobel Prize winner in economics and one of the founding fathers of AI, played a tremendous role in shaping the scientific vision and direction of the university. Simon was one of the participants in the Dartmouth workshop, where he presented the Logic Theorist, a program he had written with Allen Newell and Cliff Shaw that could solve problems in propositional logic. At Carnegie Mellon, he helped found the Graduate School of Industrial Administration, the School of Computer Science, the Department of Psychology, and the Department of Philosophy. When I arrived in 1996, logic played a prominent role in mathematics, computer science, cognitive science, philosophy, psychology, and linguistics, in no small part due to his influence. However, Carnegie Mellon is also a pioneer in machine learning, having founded the first academic department of machine learning in 2006 and the first undergraduate major in AI in 2018. On the university web pages, \href{https://ai.cmu.edu/}{AI at CMU}, you can find a history of AI at Carnegie Mellon, news and events, descriptions of various degree programs, and sample curricula. Try as you may, you will not find a single occurrence of the word ``logic.''

In 2023, the New York Times ran a series of articles designed to help the general public make sense of the revolution in AI. \href{nytimes.com/article/ai-artificial-intelligence-chatbot.html}{The first of the series} summarized the history of artificial intelligence as follows:
\begin{quote}
 A group of academics coined the term in the late 1950s as they set out to build a machine that could do anything the human brain could do---skills like reasoning, problem-solving, learning new tasks and communicating using natural language.

 Progress was relatively slow until around 2012, when a single idea shifted the entire field.

 It was called a \emph{neural network}.
\end{quote}
That brief narrative served to sum up a half-century of research in AI, carried out by some of the most talented and accomplished researchers of their time. You can't blame those of us working with symbolic methods for wondering whether our time had come.

One might object that the title of this essay is misleading because it conflates mathematics with symbolic AI. The history of mathematics consists of more than two thousand years of some of the most beautiful and creative ideas humankind has ever produced, and mathematicians generally chafe at the suggestion that formal reasoning captures the subject's essence. Mathematical thought is driven by big ideas, far-reaching intuitions, and deep insights, none of which are captured or explained by the rules of logic. With its goals of developing powerful abstractions, detecting patterns, and synthesizing disparate aspects of our experience, mathematics is as amenable to the methods of machine learning as much as it is to the methods of symbolic AI.

That said, what distinguishes mathematics among the scholarly disciplines is the level of rigor and precision with which mathematicians state and justify their claims. Although the debate over the merits of symbolic AI vs.~machine learning is not the same as the debate over the merits of mathematical vs.~scientific reasoning, both are aligned with the fundamental disagreement between Plato and Aristotle as to whether finding the right concepts is a prerequisite to knowledge or whether concepts are only what emerge from empirical data. One might argue that although mathematical ideas have supported scientific exploration and practical reasoning for centuries, their importance has now diminished. Given our limited cognitive abilities, mathematical abstraction has been a valuable tool, enabling us to squeeze out additional cognitive efficiency in our attempts to navigate a complex world. But now that we have neural networks to process the data and tell us what to make of it, the idealized mathematical representations we have been using are less helpful, pale shadows of an underlying complexity that neural networks can manage more directly. One can, therefore, argue that symbolic methods are no longer important because mathematical reasoning, as we know it, is no longer important. Technology has given us something better.

Threatened by such an argument, mathematicians may retreat to aesthetics: many of us do mathematics not because it is useful, but because we enjoy it, just as we enjoy literature and art. But aesthetic and pragmatic factors are not so easily separated. Maybe we find thinking abstractly and solving problems so enjoyable precisely because they are so generally useful; in other words, we appreciate the power of mathematics even when it is not directed at any practical goal. We might not feel the same way if mathematics were relegated to a pastime, like cricket or chess. Furthermore, an appeal to aesthetics doesn't help us solicit support for mathematics from those who do not enjoy it, and it certainly doesn't explain why we should require our children to spend so much of their time studying it.

We should take these worries seriously. I personally believe that in the age of AI, mathematics is as important as it ever has been, but that is not something we can take for granted. Contemporary AI offers us dramatically new capacities for reasoning and making decisions, and we cannot assume that methods that served us well before the advent of deep learning are still relevant today. It is therefore important to consider what role mathematics should still play in our lives, if any, and why.

\section{The Value of Mathematics}

Generally, when we ask ChatGPT a question, we care about the answer. Whenever we turn to AI for advice, we should question the wisdom of putting our lives and livelihoods in the hands of a system we don't understand. We should worry about whether the information we get is reliable, aligned with our interests, and likely to help us achieve our goals. We should worry about our safety and security, and since our decisions affect others, we should worry about whether the advice we get reflects our values and morals. For these reasons, we want AI to be transparent. We want AI to tell us not just what to do but also why; we want reasons, explanations, and justification. We expect no less from doctors, lawyers, salespeople, financial managers, and contractors, and we ought to expect the same from AI.

This highlights one role for mathematics. Mathematical language is designed to enable us to express ourselves clearly and provide rigorous justification, giving us the means to tell a system what we want and check the answers to ensure we are getting what we expect. Not everything has to be mathematized, but there are plenty of reasoning tasks where mathematical specifications and explanations are called for, and those will not go away just because we have neural networks.

This defense of mathematics, however, does not go far enough. It suggests that, in the future, the primary purpose of mathematics is for us to communicate with AI and for AI, which knows all the answers, to put them in terms we can understand. It puts AI, rather than us, at the center of the deliberative process. It fails to recognize that, when we talk about artificial intelligence, what we care about is the ability of mechanized systems to participate in \emph{our} deliberative processes and help \emph{us} reason about what we want to do and who we want to be.

Imagine you are the mayor of a town with a river running through it, and the town council has determined that, given recent growth, it is time to build a new bridge. What's a modern mayor to do? You turn to ChatGPT and say, ``Design us a bridge,'' and out come instructions for the builders. The council wants pictures, and ChatGPT obliges with images of the bridge, illuminated by the glimmer of sunrise or enveloped in a foggy winter mist. Beautiful! No mathematics is needed.

Except maybe you want blueprints for the builders, with precise specifications of the lengths and angles. They should also know the requirements the building materials must satisfy and acceptable machining tolerances. Not only do you want numeric data in the output, but you probably want to use numeric specifications in your instructions, such as the volume of traffic you want the bridge to support and how long you expect it to last. The council will undoubtedly want to know how much it will cost.

Can you trust the blueprints? The effects of a collapse would be disastrous. Until AI has established a track record, it would be wise to ask to see the calculations so that engineers can audit them and ensure they meet appropriate safety standards. They may want to check the calculations by hand or run simulations and checks using software they are more comfortable with.

However, that is not nearly enough. Building a bridge is a serious undertaking, and there is a lot to think about. What effect will the bridge have on current traffic patterns, and how will it fare with respect to anticipated growth and changes over the coming decades? Should the bridge include paths and walkways to encourage more people to walk and bike, or is it more critical to meet commercial traffic needs? How will it affect the environment, and how should we weigh environmental concerns? Will the placement of the bridge benefit some residents and harm others? What else should you take into account? These questions deserve thoughtful deliberation not just from the mayor and town council but from all the relevant stakeholders.

It's not just that we need to tell AI what we want and to make sure we get it. The point is that \emph{we often don't know what we want}, and sorting that out requires reasoning and deliberating, individually and with others. And this, in large part, is what mathematics is there for. It provides us with key capacities to reason and deliberate and to come to terms with things like measurements, costs, projections, causes and effects, likelihoods, and uncertainties.

This is a fundamental part of the human experience: thinking about the world as it is and as we want it to be, and reasoning about how to get from here to there. Being human requires thinking about what's important to us and reasoning about the consequences of our actions, and being part of society requires deliberating with others. Mathematics has evolved over centuries to support this, providing us with ways to express ourselves precisely, to produce complex chains of reasoning that we can reproduce and share with others, to evaluate our assumptions, to record and codify our reasoning processes, and to improve them over time. Mathematical concepts and abstractions matter to us because they are fundamental to the way we make sense of our experiences and the world around us. Ultimately, \emph{we} are the ones who need to decide what we should do and how we should do it. We have always enlisted technology in our efforts, but as soon as we cede our deliberative processes to technology, we will have given up something fundamental to what it means to be human.

In short, being rational means not only having goals and values but also deliberating, planning, and coordinating with others to attain them. Being able to reason about our goals and values presupposes that we can express them to ourselves and to others. For AI to help us, our interactions have to be mediated by the rich network of concepts and ideas we use to make sense of the world, and mathematics is an essential part of that network.

\section{Mathematics and AI}

I have argued in an \href{https://www.ams.org/journals/bull/2024-61-02/S0273-0979-2024-01832-1/}{article} in the \emph{Bulletin of the American Mathematical Society} that mathematics is undergoing a formal turn, as mathematicians become increasingly interested in symbolic methods and formal representations. Several high-profile formalization projects, many of them collaborative, have been carried out using the Lean Theorem Prover and its library Mathlib. The resulting digitization of mathematical content opens up new possibilities for both symbolic and neural AI.

The good news is that the value of mathematical and symbolic reasoning has not been lost on the machine-learning community. DeepMind's AlphaProof, using a reinforcement learning algorithm running on top of Lean, solved three of the six problems in this year's International Mathematics Olympiad, ranking it competitive with the world's most elite pre-college problem solvers. There has been considerable interest in improving the ability of AI to carry out mathematical reasoning. There are ongoing debates between researchers as to the extent to which symbolic systems like Lean are required for a system of AI to acquire mathematical expertise, but many feel that the most promising avenues to success require a combination of neural and symbolic methods. Even ML purists recognize that obtaining mathematical arguments, explanations, and justifications is a key goal for AI, however that can be achieved.

In recent years, therefore, the phrase ``AI for mathematics'' has come to encompass both symbolic and machine-learning methods. My colleagues and I can breathe a sigh of relief; we can also ride the AI wave, at least for now. But we won't always be able to ride it for free; we still need to show that symbolic methods have an important role to play in artificial intelligence.

Let's hope they do. Ever since Galileo proclaimed, in the seventeenth century, that the book of nature is written in the language of mathematics, mathematics has been central not only to the scientific method, but to practical decision making in several worldly pursuits, including technology, economics, finance, logistics, and public policy. Now, with the advent of AI, there are two paths we can follow. The first involves carrying out scientific reasoning and decision-making the way we have since Galileo but using AI to do it better, improving our mathematical models and obtaining a deeper understanding of their properties. The second involves bypassing mathematics, leaving AI to draw conclusions as it sees fit, and accepting its oracular conclusions. The first path opens up exciting opportunities for mathematics and science, because AI offers us new means to discover and understand phenomena that would otherwise remain opaque to us, to think and reason better, and to make better decisions. If we go down the second path, it will mean turning our back on science, relinquishing agency over our practical decisions, and giving up a vital part of what it means to be human. AI offers us the choice, but it does not tell us which path to take. It's up to us to get it right.

\paragraph{Acknowledgements.} I first grappled with these issues in a talk, with the same title as this essay, that I presented at the Santa Fe Institute in April 2023. I am grateful to the institute for giving me that opportunity, to Cris Moore for serving as an insightful respondent, and to several participants in that workshop for helpful comments and discussions. I am also grateful to Johan Commelin, Kim Morrison, and Oliver Nash for comments on a draft of this essay.

\end{document}